\documentclass[preprint,12pt,numauth]{elsarticle}

\usepackage{amssymb}
\usepackage{amsmath}

\usepackage[T1]{fontenc}
\usepackage[utf8]{inputenc}
\usepackage{array}
\usepackage{multirow}
\usepackage{amstext}
\usepackage{amsthm}
\usepackage{amsmath}
\usepackage{amssymb}
\usepackage{hyperref}
\usepackage{graphicx}
\usepackage{subcaption}
\usepackage{booktabs}
\usepackage[english]{babel}
\usepackage[table]{xcolor}

\usepackage{listings}
\usepackage{xcolor}

\definecolor{codegreen}{rgb}{0,0.6,0}
\definecolor{codegray}{rgb}{0.5,0.5,0.5}
\definecolor{codepurple}{rgb}{0.58,0,0.82}
\definecolor{backcolour}{rgb}{0.95,0.95,0.92}

\lstdefinestyle{mystyle}{
    backgroundcolor=\color{backcolour},   
    commentstyle=\color{codegreen},
    keywordstyle=\color{magenta},
    numberstyle=\tiny\color{codegray},
    stringstyle=\color{codepurple},
    basicstyle=\ttfamily\footnotesize,
    breakatwhitespace=false,         
    breaklines=true,                 
    captionpos=b,                    
    keepspaces=true,                 
    numbers=left,                    
    numbersep=5pt,                   
    showspaces=false,                
    showstringspaces=false,
    showtabs=false,                  
    tabsize=2
}

\numberwithin{equation}{section}
\numberwithin{figure}{section}
\theoremstyle{plain}
\newtheorem{thm}{\protect\theoremname}
  \theoremstyle{remark} 
  \newtheorem{remark}{Remark} 
  \theoremstyle{definition}
  
  \theoremstyle{remark}
  
  \theoremstyle{plain}
  
  \newtheorem{proposition}[thm]{\protect\propositionname}
  \providecommand{\propositionname}{Proposition}
  \newtheorem{lemma}[thm]{\protect\lemmaname}
    \providecommand{\lemmaname}{Lemma}
  
  \theoremstyle{plain}
  \newtheorem{corollary}[thm]{\protect\corollaryname}
  \providecommand{\corollaryname}{Corollary}

\makeatother

\usepackage{babel}
  \providecommand{\conjecturename}{Conjecture}
  \providecommand{\definitionname}{Definition}
  \providecommand{\examplename}{Example}
  \providecommand{\remarkname}{Remark}
\providecommand{\lemmaname}{Lemma}
\providecommand{\theoremname}{Theorem}

\journal{Number theory}

\begin{document}

\begin{frontmatter}

\title{On the Stieltjes Approximation Error to the Logarithmic Integral} 

\author{Jonatan Gomez} 

\affiliation{organization={Department of Computer Engineering, Universidad Nacional de Colombia},
            addressline={Carrera 30 No 45-03}, 
            city={Bogota},
            postcode={110911}, 
            state={DC},
            country={Colombia}}

\begin{abstract}
We study the approximation errors $\varepsilon(x)=\operatorname{li}_{*}(x)-\operatorname{li}(x)$ and $\varepsilon_{S}(x)=\operatorname{li}_{S}(x)-\operatorname{li}(x)$ arising from the classical Stieltjes--van Boven and pure Stieltjes asymptotic expansions of the logarithmic integral, respectively. Our analysis establishes a top-down reduction based on the discrete error sequence $\varepsilon_n=\varepsilon(e^n)$ and its localized partial variations, for which we derive new, unconditional analytic bounds. By developing precise integral representations for these error increments and establishing sharp, fully explicit upper and lower estima-tes for the underlying kernel integrals, we obtain uniform bounds valid for every integer $n\geq 1$. These discrete results are subsequently lifted to the real domain to provide comprehensive estimates for all real $x\geq e$. In particular, we prove that $\varepsilon(x) = \frac{1}{3} \sqrt{2\pi/\log x} + O(1/\log x)$ as $x \to \infty$. Furthermore, we deploy these continuous bounds to uncover the fine structural behavior of the pure Stieltjes expansion, demonstrating that $|\varepsilon_{S}(x)| = \sqrt{2\pi}\,|1/3-\alpha|/\sqrt{\log x} + O(1/\log x)$ as $x\to\infty$, where $n=\lfloor\log x\rfloor$ and $\alpha=\log x-n$. This establishes a novel and explicit asymptotic description of the absolute pure Stieltjes error, fully capturing its non-trivial oscillatory dependence on the fractional part of the logarithm.
\end{abstract}



\end{frontmatter}

\section{Introduction and Historical Background} \label{sec:intro}

    The logarithmic integral, denoted by $\operatorname{li}(x)$, is a classical special function in mathematics with applications in physics and number theory. In particular, it plays a central role in number theory as an approximation to the prime-counting function $\pi(x)$, which counts the number of primes less than or equal to a given real number $x$~\cite{narkiewicz2000}. For $x>1$, the logarithmic integral is defined by
    
    \[ \operatorname{li}(x) = \operatorname{PV}\int_0^x \frac{dt}{\log t}, \]

    where $\operatorname{PV}$ denotes the Cauchy principal value at $t=1$.

    As $x\to\infty$, the logarithmic integral admits the classical asymptotic expansion in inverse powers of $\log x$,

    \[ \operatorname{li}(x) \sim \frac{x}{\log x} \sum_{k=0}^{\infty}\frac{k!}{\log^k x}. \]

    This asymptotic expansion is divergent and must therefore be truncated after finitely many terms. In particular, for any fixed positive integer $m$, its $m$-term truncation is
    
    \[ \operatorname{li}_m(x) = \frac{x}{\log x} \sum_{k=0}^{m-1}\frac{k!}{\log^k x}, \]

    and satisfies~\cite{masina2019useful}

    \[ \operatorname{li}(x) = \operatorname{li}_m(x) + O\left(\frac{m!x}{\log^{m+1}x}\right), \qquad x\to\infty. \]

    Stieltjes~\cite{ASENS_1886_3_3__201_0} observed that the truncation index minimizing the approximation error is asymptotically given by the optimal choice $n^*=\lfloor\log x\rfloor$. Accordingly, the pure Stieltjes approximation to the logarithmic integral is defined, for $x\geq e$, by

    \[ \operatorname{li}_{S}(x) = \frac{x}{\log x} \sum_{k=0}^{n^*-1}\frac{k!}{\log^k x}, \qquad n^*=\lfloor\log x\rfloor, \]

    with an associated approximation error denoted by $\varepsilon_{S}(x) = \operatorname{li}_{S}(x)-\operatorname{li}(x)$.

    Van Boven~\cite{Boven2012AsymptoticSO} introduced a refinement of this classical expansion by incorporating a fractional remainder term. The resulting Stieltjes--van Boven (SV) approximation is defined, for $x\geq e$, as

    \[ \operatorname{li}_{*}(x) = \frac{(\log x-n^*)x n^*!}{\log^{n^*+1}x} + \frac{x}{\log x} \sum_{k=0}^{n^*-1}\frac{k!}{\log^k x}, \qquad n^*=\lfloor\log x\rfloor. \]

    We denote its corresponding approximation error by $\varepsilon(x) = \operatorname{li}_{*}(x)-\operatorname{li}(x)$. 

    Although Stieltjes~\cite{ASENS_1886_3_3__201_0} originally suggested that the order of approximation for these errors should be proportional to $\sqrt{2\pi/\log x}$, subsequent observations by van Boven~\cite{Boven2012AsymptoticSO} indicated that this classical estimate might overstate the true magnitude of the error. A rigorous, unconditional description of the precise asymptotic behavior and qualitative properties of both $\varepsilon(x)$ and $\varepsilon_0(x)$ has remained an open problem.

    The present work provides a complete and unconditional solution to this problem. First, we establish the strict positivity and monotonicity properties of the Stieltjes--van Boven approximation error $\varepsilon(x)$ and prove the asymptotic formula conjectured in earlier versions of this work:
    
    \[ \varepsilon(x) = \frac{1}{3}\sqrt{\frac{2\pi}{\log x}} + O\left(\frac{1}{\log x}\right), \qquad x\to\infty. \]
    
    Second, we deploy these sharp continuous estimates to uncover the fine structure of the absolute error associated with the pure Stieltjes expansion. Writing $n=\lfloor\log x\rfloor$ and letting $\alpha=\log x-n$ represent the fractional part of the logarithm, we prove that

    \[ |\varepsilon_{S}(x)| = \sqrt{\frac{2\pi}{\log x}} \left|\frac{1}{3}-\alpha\right| + O\left(\frac{1}{\log x}\right), \qquad x\to\infty. \]

    This result establishes a novel, explicit asymptotic description that captu-res how the absolute error under the pure Stieltjes expansion undergoes non-trivial oscillations driven entirely by the fractional remainder $\alpha$.

    This paper represents a comprehensive reorganization and expansion of our previous manuscript (\href{https://arxiv.org/abs/2406.12152v2}{arXiv:2406.12152v2}). Beyond strengthening the underlying error bounds to make them fully unconditional, the presentation has been restructured into a modular, top-down architecture to separate the global analytical results from the extensive machinery required to prove them.

    The remainder of this paper is organized as follows. Section~\ref{sec:main_results} establishes our framework and formally states the main theorems, including explicit bounds valid for all $x \ge e$. Section~\ref{sec:discrete_reduction} introduces the top-down strategy by discretizing the error functions at $x = e^k$ and deriving exact integral representations for their increments. The necessary technical machinery is then compartmentalized: Section~\ref{sec:kernel_estimates} deduces sharp estimates for the underlying kernel integrals, while Section~\ref{sec:factorial_estimates} provides refinements on factorial factors and Stirling-type expansions. Section~\ref{sec:proofs} synthesizes these components to complete the proofs of the main theorems.
\section{Main Results} \label{sec:main_results}

In this section, we formally present the main contributions of this work. We establish both the global asymptotic behaviors as $x \to \infty$ and sharp, unconditional explicit bounds. Throughout this paper, for any real number $x \ge e$, we denote the optimal truncation index by $n = \lfloor\log x\rfloor$ and its corresponding fractional part by $\alpha = \log x - n$, so that $\alpha \in [0, 1)$.

Our analysis initiates in the discrete domain by evaluating the errors at the exponential lattice points $x = e^n$. At these specific points, the fractional remainder vanishes, causing the Stieltjes--van Boven error and the pure Stieltjes error to coincide. We define this discrete sequence for $n \in \mathbb{N}_{\ge 1}$ as:
\begin{equation} \label{eq:discrete_coincidence}
    \varepsilon_n = \varepsilon(e^n) = \varepsilon_S(e^n).
\end{equation}

Our first main result provides sharp, non-trivial, and fully unconditional explicit bounds for this fundamental discrete sequence.

\begin{thm}[Sharp Explicit Bounds for the Discrete Error Sequence] \label{thm:varepsilon_n}
For every positive integer $n \ge 1$, the discrete error value $\varepsilon_n$ satisfies the explicit inequalities:
\begin{equation} \label{eq:discrete_bounds_formula}
     0 \le \frac{\sqrt{2\pi}}{3\sqrt{n}}  - \frac{\sqrt{2\pi}\left(\sqrt{n+1}-1\right)}{3\sqrt{n}\sqrt{n+1}\left(\sqrt{n}+\sqrt{n+1}\right)} \le \varepsilon_{n} \le \frac{\sqrt{2\pi}}{3\sqrt{n}}.
\end{equation}
\end{thm}

By leveraging the precise discrete control established in Theorem \ref{thm:varepsilon_n} alongside the analytical behavior of the localized error functions between lattice points, we lift these properties to the continuous real domain. Our second main theorem establishes the global asymptotic behavior and qualitative properties of the Stieltjes--van Boven approximation error $\varepsilon(x)$ for all real arguments.

\begin{thm}[Asymptotic Behavior and Properties of the SV Error] \label{thm:sv_main}
For every real number $x \ge e$, the Stieltjes--van Boven approximation error $\varepsilon(x) = \operatorname{li}_{*}(x) - \operatorname{li}(x)$ is strictly positive and monotonically decreasing. Furthermore, as a consequence of Theorem \ref{thm:varepsilon_n}, it satisfies the asymptotic expansion:
\begin{equation} \label{eq:sv_asymptotic}
    \varepsilon(x) = \frac{1}{3}\sqrt{\frac{2\pi}{\log x}} + O\left(\frac{1}{\log x}\right), \qquad x \to \infty.
\end{equation}
\end{thm}

Crucially, the continuous control of $\varepsilon(x)$ allows us to unearth the fine structural behavior of the pure Stieltjes error $\varepsilon_S(x)$. Unlike the smooth decay of the SV error, the pure Stieltjes error exhibits asymptotic oscillations driven non-trivially by the fractional remainder $\alpha$.

\begin{thm}[Fine Structure and Oscillations of the Pure Stieltjes Error] \label{thm:stieltjes_main}
Let $x \ge e$, $n = \lfloor\log x\rfloor$, and $\alpha = \log x - n$. The absolute approximation error associated with the pure Stieltjes expansion, $\varepsilon_S(x) = \operatorname{li}_S(x) - \operatorname{li}(x)$, satisfies:
\begin{equation} \label{eq:stieltjes_asymptotic}
    |\varepsilon_{S}(x)| = \sqrt{\frac{2\pi}{\log x}} \left|\frac{1}{3}-\alpha\right| + O\left(\frac{1}{\log x}\right), \qquad x \to \infty.
\end{equation}
\end{thm}

\begin{remark}
Theorem \ref{thm:stieltjes_main} confirms that the classical error order $\sqrt{2\pi/\log x}$ originally suggested by Stieltjes~\cite{ASENS_1886_3_3__201_0} is structurally modulated by the term $|1/3 - \alpha|$. This implies that the pure Stieltjes error shrinks significantly near the points where $\alpha \approx 1/3$, explaining why the classical unmodulated bound overstates the actual error as observed by van Boven~\cite{Boven2012AsymptoticSO}.
\end{remark}

To establish the bridge between these continuous theorems and the analytical discrete machinery developed in the subsequent sections, we conclude by highlighting the strategic roadmap of our proofs. 

The explicit bounds in Theorem \ref{thm:varepsilon_n} are achieved via a top-down reduction strategy. In Section \ref{sec:discrete_reduction}, we introduce the successive discrete decrements $\Delta_{n+1} = \varepsilon_n - \varepsilon_{n+1}$ and derive their exact integral representations. The underlying technical proof is then compartmentalized: Section \ref{sec:kernel_estimates} deduces sharp upper and lower estimates for the associated kernel integrals, while Section \ref{sec:factorial_estimates} provides localized refinements on factorial factors and Stirling-type expansions. Section \ref{sec:proofs} synthesizes these components to complete the proofs of Theorem \ref{thm:varepsilon_n}, Theorem \ref{sec:main_results}, and Theorem \ref{thm:stieltjes_main}.

\section{The Partial Error Function and Integral Representations} \label{sec:discrete_reduction}

To establish the sharp continuous estimates asserted in the main theorems, we must bridge the behavior of the error function between the discrete exponential lattice points. We achieve this by introducing a localized, continuous partial approximation error function, denoted as $\Delta(x)$. The following proposition establishes the exact structural link that allows us to decompose the continuous error into its discrete baseline and a localized integral correction.

\begin{proposition} \label{prop:varepsilon(x)-delta(x)}
For any real number $x > e$, the Stieltjes--van Boven continuous approximation error satisfies the decomposition:
\begin{equation} \label{eq:error_decomposition_identity}
    \varepsilon(x) = \varepsilon_n -\Delta(x),
\end{equation}  
where $n = \lceil \log(x) \rceil - 1$ and the partial error function $\Delta(x)$ is defined by:
\begin{equation} \label{eq:partial_error_def}
    \Delta(x)=\int_{e^{n}}^{x}\frac{n!}{\log^{n+1} t}\,dt - \frac{(\log(x)-n)x n!}{\log^{n+1}(x)}.
\end{equation}  
\end{proposition}

\begin{proof}
Let $n = \lceil \log(x) \rceil - 1$ for any real number $x > e$. Notice that $n=\lfloor\log(x)\rfloor$ if $x\ne e^k$ for all $k \in \mathbb{N}_{\ge 1}$, and $n=k-1$ if $x = e^k$ for some positive integer $k$. By the definition of the logarithmic integral and the additivity of definite integrals, we split the domain as:
\begin{equation*}
    \operatorname{li}(x) = \operatorname{li}(e^n) + \int_{e^n}^{x}\frac{dt}{\log t}.
\end{equation*}
Applying integration by parts $n$ times to the remaining integral yields:
\begin{equation*}
    \operatorname{li}(x) = \operatorname{li}(e^{n}) + \frac{x}{\log(x)}\sum_{k=0}^{n-1}\frac{k!}{\log^{k}(x)} - \frac{e^{n}}{n}\sum_{k=0}^{n-1}\frac{k!}{n^{k}} + \int_{e^{n}}^{x}\frac{n!}{\log^{n+1} t}\,dt.
\end{equation*}
By adding and subtracting the fractional remainder term $\frac{(\log(x)-n)x n!}{\log^{n+1}(x)}$, and identifying the structural definition of the Stieltjes--van Boven approximation $\operatorname{li}_{*}(x)$, we obtain:
\begin{equation*}
    \operatorname{li}(x) = \operatorname{li}(e^{n}) + \operatorname{li}_{*}(x)  - \operatorname{li}_{*}(e^n) + \int_{e^{n}}^{x}\frac{n!}{\log^{n+1} t}\,dt -  \frac{(\log(x)-n)x n!}{\log^{n+1}(x)}.
\end{equation*}
Rearranging the terms allows us to collect the functional differences:
\begin{equation*}
    \operatorname{li}_{*}(x) - \operatorname{li}(x) = \operatorname{li}_{*}(e^{n}) -\operatorname{li}(e^{n}) - \left[\int_{e^{n}}^{x}\frac{n!}{\log^{n+1} t}\,dt - \frac{(\log(x)-n)x n!}{\log^{n+1}(x)}\right].
\end{equation*}
By utilizing the definition of the continuous error function $\varepsilon(x) = \operatorname{li}_{*}(x) - \operatorname{li}(x)$ and noting that at the lattice boundary $\varepsilon(e^n) = \varepsilon_n$, this identity reduces to:
\begin{equation*}
    \varepsilon(x) = \varepsilon_{n} - \left[\int_{e^{n}}^{x}\frac{n!}{\log^{n+1} t}\,dt - \frac{(\log(x)-n)x n!}{\log^{n+1}(x)}\right].
\end{equation*}
Therefore, substituting the explicit definition of $\Delta(x)$ given in \eqref{eq:partial_error_def} yields the targeted relation:
\begin{equation*}
    \varepsilon(x) = \varepsilon_{n} - \Delta(x),
\end{equation*}
which completes the proof of the proposition.
\end{proof}

Our next objective is to derive an exact, tractable integral representation for $\Delta(x)$ that eliminates the singularities and reveals its underlying structural kernel.

\begin{lemma} \label{lem:Delta(x)inu}
For all $x > e$ and $n = \lceil \log(x) \rceil - 1$, the partial error function $\Delta(x)$ admits the integral representation:
\begin{equation} \label{eq:compact_integral_form}
    \Delta(x) = \int_{0}^{\log(x)-n} f_n(u) \,du, \quad \text{where} \quad f_n(u) =\frac{u(1-u)\,e^{n+u}\,n!}{(n+u)^{\,n+2}}.
\end{equation}
\end{lemma}

\begin{proof}
Differentiating the second term of \eqref{eq:partial_error_def} with respect to $x$ yields:
\begin{equation*}
    \left(\frac{(\log(x)-n)\, x\, n!}{\log^{n+1}(x)}\right)' = \frac{n!}{\log^{n+1}(x)} \left[ 1+(\log(x)-n) -\frac{(n+1)(\log(x)-n)}{\log(x)} \right].
\end{equation*}
Notice that this second term in $\Delta(x)$ can be fully recovered by computing the integral of its derivative over the interval $[e^n, x]$:
\begin{equation*}
    \int_{e^n}^{x} \frac{n!}{\log^{n+1}(t)} \left( 1+(\log(t)-n) -\frac{(n+1)(\log(t)-n)}{\log(t)} \right)\,dt.
\end{equation*}
Incorporating this expression back into the definition of $\Delta(x)$ in \eqref{eq:partial_error_def}, we obtain:
\begin{equation*}
    \Delta(x) = \int_{e^{n}}^{x} - \frac{n!}{\log^{n+1}(t)} \left( (\log(t)-n) -\frac{(n+1)(\log(t)-n)}{\log(t)} \right)\,dt.
\end{equation*}
Simplifying and reorganizing the integrand leads directly to:
\begin{equation*}
    \Delta(x) = \int_{e^{n}}^{x} \frac{n!}{\log^{n+2}(t)} (\log(t)-n)(\log(t)-n-1)\,dt.
\end{equation*}
Finally, introducing the change of variables $u=\log(t)-n$, so that $dt = e^{n+u}du$, yields the compact form:
\begin{equation*}
    \Delta(x) = \int_{0}^{\log(x)-n} f_n(u) \,du, \quad \text{with }\quad f_n(u) =\frac{u(1-u)\,e^{n+u}\,n!}{(n+u)^{\,n+2}}
\end{equation*}
which holds for all $x > e$ and $n = \lceil \log(x) \rceil - 1$, completing the proof.
\end{proof}

By analyzing the sign and structural properties of the integrand $f_n(u)$, we readily establish the monotonicity of the partial error function and its precise connection to the discrete boundary decrements $\Delta_{n+1}$.

\begin{proposition} \label{prop:Delta(x)order}
For all integers $n\ge1$ and all real numbers $e^{n} < x < y \le e^{n+1}$, the partial error function satisfies:
\begin{equation} \label{eq:monotonicity_chain}
0 <\Delta(x) < \Delta(y) \le \Delta(e^{n+1})=\Delta_{n+1}.
\end{equation}
\end{proposition}
    
\begin{proof}
Since the kernel function satisfies $f_n(u)>0$ for all $u\in(0,1)$ and vanishes at the boundaries $f_n(0)=f_n(1)=0$, the integral \eqref{eq:compact_integral_form} is strictly positive and strictly increasing for every $e^n < x \le e^{n+1}$. This directly implies the monotonicity of $\Delta(x)$. Consequently, evaluating the function at the upper boundary $x = e^{n+1}$ yields the maximum value, which matches the discrete decrement:
\begin{equation*}
0 <\Delta(x) < \Delta(y) \le \Delta(e^{n+1})=\Delta_{n+1}
\end{equation*}
for all $n\ge1$ and all real numbers $e^{n} < x < y \le e^{n+1}$.
\end{proof}

The monotonic structure established in Proposition \ref{prop:Delta(x)order} implies that bounding the continuous function $\Delta(x)$ reduces to controlling the discrete endpoint decrement $\Delta_{n+1} = \int_0^1 f_n(u)du$. To achieve this without losing precision, the integrand $f_n(u)$ must be systematically decoupled. As outlined in the introduction, the sharp estimates for the kernel integrals and the factorial modulations will be individually structured in Section \ref{sec:kernel_estimates} and Section \ref{sec:factorial_estimates}, respectively, before establishing the definitive bounds for $\Delta_{n+1}$.

\section{Analysis and Estimates of Kernel Integrals} \label{sec:kernel_estimates}

This section establishes sharp, uniform algebraic bounds for the core transcendental integrals that govern the behavior of the error. As established in Proposition \ref{prop:Delta(x)order}, the evaluation of the error over the exponential lattice points depends on the endpoint decrements $\Delta_{n+1} = \int_{0}^{1} f_n(\alpha) \, d\alpha$. 

After isolating the factorial components via Stirling-type modulations (which will be rigorously detailed in Section \ref{sec:factorial_estimates}), we are led to the analysis of two fundamental kernel integrals, denoted for $n \ge 10$ as:
\begin{equation} \label{eq:Jn_definition}
    J_n = \int_0^1 \frac{ e^{\frac{1}{12n+1}-\frac{1}{16n} - \frac{\alpha(1-\alpha)}{4n}} \alpha(1-\alpha)}{(n+\alpha)^{3/2}}\,d\alpha,
\end{equation}
and
\begin{equation} \label{eq:Ln_definition}
    L_n = \int_0^1 \frac{e^{\frac{1}{12n+1} - \frac{\alpha(1-\alpha)}{4n}}\alpha(1-\alpha)}{(n+\alpha)^{3/2}}\,d\alpha.
\end{equation}

The following lemmas reduce these transcendental integrals to explicit algebraic expressions in terms of the baseline function $B_n$ and the secondary error residuals $S_n$ and $R_n$.

\begin{lemma}[Kernel Lower Bound] \label{lem:integral-lower}
For every integer $n \ge 10$, the kernel integral $J_n$ satisfies:
\begin{equation} \label{eq:Jn_lower_bound}
    B_n - S_n \le J_n,
\end{equation}
where
\begin{equation} \label{eq:Bn_Sn_def}
    B_n = \frac{1}{3\sqrt{n}(n+1)\left(1+\sqrt{\frac{n}{n+1}}\right)} \quad \text{and} \quad S_n = \frac{1}{n^{5/2}} \left( \frac{7}{1440}+\frac{1}{8n} \right).
\end{equation}
\end{lemma}

\begin{proof}
Let $x=\frac{1}{n}$. For $0\le\alpha\le1$ and $n \ge 10$, the elementary inequalities $(1+x\alpha)^{-3/2} \ge 1-\frac{3}{2}x\alpha$ and $e^y\ge 1+y$ allow us to bound the integrand of \eqref{eq:Jn_definition} from below as:
\begin{equation*}
    J_n\ge \frac{1}{n^{3/2}}\int_0^1 \alpha(1-\alpha) \left(1-\frac{3\alpha}{2n}\right) \left( 1 + \left(\frac{1}{12n+1}-\frac{1}{16n}\right) - \frac{\alpha(1-\alpha)}{4n}\right)\,d\alpha.
\end{equation*}
Both factors on the right-hand side are non-negative for $n\ge 10$, meaning that expanding the product preserves the lower bound. Letting $A_n=\frac{1}{12n+1}-\frac{1}{16n}$ and evaluating the explicit integrals $\int_0^1\alpha(1-\alpha)\,d\alpha=\frac{1}{6}$, $\int_0^1\alpha^2(1-\alpha)\,d\alpha=\frac{1}{12}$, $\int_0^1\alpha^2(1-\alpha)^2\,d\alpha=\frac{1}{30}$, and $\int_0^1\alpha^3(1-\alpha)^2\,d\alpha=\frac{1}{60}$, we obtain:
\begin{equation*}
    J_n\ge\frac{1}{n^{3/2}}\bigg( \frac{1}{6}+\frac{A_n}{6}-\frac{1}{120n} - \frac{1+A_n}{8n}+\frac{1}{160n^2} \bigg).
\end{equation*}
By writing $A_n=\frac{1}{48n}-\frac{1}{12n(12n+1)}$, it follows that $\frac{1}{48n}-\frac{1}{144n^2} \le A_n\le\frac{1}{48n}$. Applying the lower bound for the positive term $A_n/6$ and the upper bound for the negative term $-A_n/(8n)$, we find:
\begin{equation*}
\begin{aligned}
    J_n &\ge \frac{1}{n^{3/2}}\left[ \frac{1}{6} + \frac{1}{288n}-\frac{1}{864n^2} - \frac{1}{120n} - \frac{1}{8n} - \frac{1}{384n^2} + \frac{1}{160n^2} \right] \\ &\ge \frac{1}{n^{3/2}}\left( \frac{1}{6}-\frac{187}{1440n} + \frac{43}{17280n^2} \right).    
\end{aligned}
\end{equation*}

Discarding the positive $x^2$ term yields the clean lower estimate:

\begin{equation} \label{eq:Jn_intermediate_1}
    J_n \ge \frac{1}{6n^{3/2}} - \frac{187}{1440n^{5/2}}.
\end{equation}

Alternatively, the baseline expression $B_n$ can be expressed as 

\[B_n = \frac{1}{3n^{3/2} \left(1+x+\sqrt{1+x}\right)}.\] 

For $0\le x\le1$, utilizing the Taylor expansion of the square root yields $\sqrt{1+x}\ge 1+\frac{x}{2}-\frac{x^2}{8}$. Since $3\left(\frac{1}{6}-\frac{x}{8}+\frac{x^2}{8}\right) \left(2+\frac{3x}{2}-\frac{x^2}{8}\right)-1 = \frac{x^2(8+39x-3x^2)}{64} \ge 0$, we establish the algebraic inequality:
\begin{equation*}
    \frac{1}{3(1+x+\sqrt{1+x})} \le \frac{1}{6}-\frac{x}{8}+\frac{x^2}{8}.
\end{equation*}

Setting $x=1/n$ gives the following upper bound for our target baseline:

\begin{equation} 
    \label{eq:Bn_intermediate_2}
    B_n\le\frac{1}{n^{3/2}} \left(\frac{1}{6}-\frac{1}{8n}+\frac{1}{8n^2}\right).
\end{equation}

Finally, combining \eqref{eq:Jn_intermediate_1} and \eqref{eq:Bn_intermediate_2}, we isolate the difference:

\begin{equation*}
    B_n-J_n \le \frac{1}{n^{5/2}} \left( \frac{187}{1440}-\frac{1}{8}+\frac{1}{8n} \right) = \frac{1}{n^{5/2}} \left( \frac{7}{1440}+\frac{1}{8n} \right) = S_n.
\end{equation*}

Rearranging terms yields $B_n - S_n \le J_n$, completing the proof.
\end{proof}

\begin{lemma}[Kernel Upper Bound] \label{lem:integral-upper}
For every integer $n\ge 10$, the kernel integral $L_n$ satisfies:
\begin{equation} \label{eq:Ln_upper_bound}
    L_n \le B_n + R_n,
\end{equation}
where
\begin{equation} \label{eq:Rn_def}
    R_n = \frac{1}{n^{5/2}}\left(\frac{49}{720}+\frac{97}{1000n}\right).
\end{equation}
\end{lemma}

\begin{proof}
Let $x=\frac{1}{n}$ and $w=\alpha(1-\alpha)$. We rewrite \eqref{eq:Ln_definition} compactly as $L_n=\frac{1}{n^{3/2}}\int_0^1 w\,e^z(1+x\alpha)^{-3/2}\,d\alpha$, where $z=\frac{x}{12+x}-\frac{xw}{4}$. Since $0\le w\le1/4$ for $\alpha \in [0,1]$, the exponent satisfies $0\le z\le\frac{x}{12}$. Applying Taylor's theorem to $e^z$ leads to:
\begin{equation*}
    e^z\le 1+z+\frac{11}{1000}x^2 \le 1+x\left(\frac{1}{12}-\frac{w}{4}\right) + \frac{11}{1000}x^2.
\end{equation*}
Similarly, utilizing the binomial expansion for the algebraic core, we have $(1+y)^{-3/2} \le 1-\frac{3}{2}y+\frac{15}{8}y^2$ for $y \ge 0$, which implies:
\begin{equation*}
    (1+x\alpha)^{-3/2} \le 1-\frac{3}{2}x\alpha+\frac{15}{8}x^2\alpha^2.
\end{equation*}

Setting $b=\frac{1}{12}-\frac{w}{4}$ and multiplying these two uniform bounds, we obtain the joint expansion:

\begin{equation*}
\begin{aligned}
e^z(1+x\alpha)^{-3/2}
\le {}&
1+x\left(b-\frac32\alpha\right)\\
&+x^2\left(
\frac{11}{1000}
+\frac{15}{8}\alpha^2
-\frac32b\alpha
\right)\\
&+x^3\left(
\frac{15}{8}b\alpha^2
-\frac{33}{2000}\alpha
\right)
+\frac{165}{8000}\alpha^2x^4 \\ \le {}& 1+x\left(b-\frac{3}{2}\alpha\right) +x^2\left( \frac{11}{1000} + \frac{15}{8}\alpha^2 - \frac{3}{2}b\alpha \right) \\ &+ \frac{15}{8}b\alpha^2x^3 + \frac{165}{8000}\alpha^2x^4 - \frac{33}{2000}\alpha x^3.
\end{aligned}
\end{equation*}

Since $0\leq\alpha\leq 1$, $b=\frac1{12}-\frac{\alpha(1-\alpha)}4
\geq \frac1{48}>0$, and $x>0$, both $-\frac{3}{2}b\alpha x^2$ and $-\frac{33}{2000}\alpha x^3$ are non-positive. Hence, both terms may be discarded when deriving an upper bound. For $n\ge10$, we have $0<x\le1/10$, allowing us to bound the higher-order error terms as:
\begin{equation*}
    \frac{15}{8}b\alpha^2x^3 + \frac{165}{8000}\alpha^2x^4 \le \left( \frac{15}{8}\cdot\frac{1}{12}\cdot\frac{1}{10} + \frac{165}{8000}\cdot\frac{1}{100} \right)\alpha^2x^2 < \frac{21}{1000}\alpha^2x^2.
\end{equation*}

Gathering the coefficients simplifies the inequality to:
\begin{equation*}
    e^z(1+x\alpha)^{-3/2} \le 1+x\left(\frac{1}{12}-\frac{w}{4}-\frac{3}{2}\alpha\right) + x^2\left(\frac{11}{1000}+\frac{1896}{1000}\alpha^2\right).
\end{equation*}
Multiplying this uniform bound by $w$ and integrating term by term using the standard moments $\int_0^1w\,d\alpha=\frac{1}{6}$, $\int_0^1w^2\,d\alpha=\frac{1}{30}$, $\int_0^1\alpha w\,d\alpha=\frac{1}{12}$, and $\int_0^1\alpha^2w\,d\alpha=\frac{1}{20}$, we find:
\begin{equation} \label{eq:Ln_intermediate_final}
\begin{aligned}
    L_n &\le \frac{1}{n^{3/2}} \left[ \frac{1}{6}-\frac{43}{360n} + \frac{1}{n^2}\left(\frac{11}{6000}+\frac{1896}{20000}\right)\right] \\ &\le \frac{1}{n^{3/2}}\left(\frac{1}{6}-\frac{41}{720n} + \frac{97}{1000n^2}\right).
\end{aligned} 
\end{equation}

Conversely, applying $\sqrt{1+x}\le1+x/2$ and $(1+y)^{-1}\ge1-y$ for $y\ge0$ directly to the explicit definition of $B_n$ yields the lower baseline estimation:

\begin{equation} \label{eq:Bn_lower_algebraic}
\begin{aligned}
    B_n &\ge \frac{1}{3n^{3/2}(1+x+1+x/2)} \\ &\ge \frac{1}{6n^{3/2}}\left(1-\frac{3x}{4}\right) \\ &\ge \frac{1}{n^{3/2}} \left(\frac{1}{6}-\frac{1}{8n}\right).
\end{aligned}
\end{equation}
Subtracting \eqref{eq:Bn_lower_algebraic} from \eqref{eq:Ln_intermediate_final} isolates the total upper residual:
\begin{equation*}
    L_n-B_n \le \frac{1}{n^{5/2}} \left(\frac{1}{8}-\frac{41}{720} + \frac{97}{1000n}\right) = \frac{1}{n^{5/2}}\left(\frac{49}{720}+\frac{97}{1000n}\right) = R_n,
\end{equation*}
which completes the proof.
\end{proof}

With the uniform control over the localized factorial ratios firmly established in Proposition \ref{prop:n_(n+a)}, we are now in a position to synthesize our preliminary estimates. The following lemma provides the ultimate analytical closure for the discrete endpoint decrements $\Delta_{n+1}$. By merging the combinatorial bounds of this section with the transcendental kernel integral estimates from Section \ref{sec:kernel_estimates}, we reduce the behavior of these discrete steps to purely explicit algebraic functions.

\begin{lemma} \label{lem:Delta_n}
For all integers $n \ge 10$, the discrete endpoint decrements $\Delta_{n+1}$ satisfy the uniform bounds:
\begin{equation} \label{eq:Delta_n_final_bounds}
   \sqrt{2\pi}\left(B_{n} - S_{n}\right) \le \Delta_{n+1} \le  \sqrt{2\pi}\left(B_{n} + R_{n}\right), 
\end{equation}
where the algebraic functions $B_n$, $S_n$, and $R_n$ are explicitly defined as in \eqref{eq:Bn_Sn_def} and \eqref{eq:Rn_def}.
\end{lemma}

\begin{proof}
Consider the continuous localized kernel function $f_n(\alpha)$ introduced in Lemma \ref{lem:Delta(x)inu}. For $n \ge 10$ and $\alpha \in [0,1]$, we can rearrange its components to isolate the shifting factorial core as follows:
\begin{equation*}
    f_n(\alpha) =\frac{\alpha(1-\alpha)\,e^{n+\alpha}\,n!}{(n+\alpha)^{\,n+2}} = \frac{\alpha(1-\alpha)}{n+\alpha}\cdot\frac{e^{n+\alpha}\,n!}{(n+\alpha)^{n+1}}.
\end{equation*}
Invoking the sharp uniform bounds for the localized factorial ratio established in Proposition \ref{prop:n_(n+a)}, we can substitute the transcendental factorial factor to obtain the direct integrand inequalities:
\begin{equation*}
    \frac{\sqrt{2\pi}e^{\frac{1}{12n+1}-\frac{\alpha(1-\alpha)}{4n} -\frac{1}{16n}}\alpha(1-\alpha)}{(n+\alpha)^{\frac{3}{2}}} \le f_n(\alpha) \le \frac{\sqrt{2\pi}e^{\frac{1}{12n}-\frac{\alpha(1-\alpha)}{4n}}\alpha(1-\alpha)}{(n+\alpha)^{\frac{3}{2}}}.
\end{equation*}
Integrating this uniform chain over the unit interval $[0,1]$ preserves the direction of the inequalities, yielding:
\begin{equation*}
    \sqrt{2\pi} J_n  \le \int_0^{1} f_n(\alpha)\, d\alpha \le  \sqrt{2\pi} L_n, 
\end{equation*}
where the transcendental kernel integrals $J_n$ and $L_n$ match exactly definitions \eqref{eq:Jn_definition} and \eqref{eq:Ln_definition}:
\begin{equation*}
    J_n = \int_0^{1}\frac{e^{\frac{1}{12n+1}-\frac{\alpha(1-\alpha)}{4n} - \frac{1}{16n}}\alpha(1-\alpha)}{(n+\alpha)^{\frac{3}{2}}} \, d\alpha \quad \text{and} \quad L_n = \int_0^{1} \frac{e^{\frac{1}{12n}}\alpha(1-\alpha)}{(n+\alpha)^{\frac{3}{2}}}\, d\alpha.
\end{equation*}
By Proposition \ref{prop:Delta(x)order}, the central integral corresponds precisely to the discrete endpoint decrement, meaning $\Delta_{n+1} = \int_0^1 f_n(\alpha)d\alpha$. Finally, substituting the explicit algebraic evaluations for $J_n$ and $L_n$ proven in Lemma \ref{lem:integral-lower} and Lemma \ref{lem:integral-upper}, we arrive at the targeted unconditional bounds:
\begin{equation*}
    \sqrt{2\pi}\left(B_{n} - S_{n}\right) \le \Delta_{n+1} \le \sqrt{2\pi}\left(B_{n} + R_{n}\right),
\end{equation*}
which holds for all integers $n \ge 10$, completing the proof.
\end{proof}

\section{Factorial Estimates and Stirling Refinements} \label{sec:factorial_estimates}

To achieve fully unconditional and sharp explicit bounds for the Stieltjes expansion error, it is necessary to establish precise localized estimates for the factorial function. In the integral representation of the decrements developed in Section \ref{sec:discrete_reduction}, the core factorial ratio appears under the following structural form:
\begin{equation} \label{eq:factorial_ratio_target}
    \frac{e^{n+\alpha}\,n!}{(n+\alpha)^{n+1}}.
\end{equation}
The purpose of this section is to isolate this term and derive uniform, non-asymptotic algebraic bounds valid for all real numbers $\alpha \in [0,1]$ and integers $n \ge 1$. This allows us to master the discrete shifting corrections and boundary effects.

Our analysis builds upon elementary logarithm series expansions and the analytical properties of a localized exponential difference function, which we formalize through the subsequent lemmas.

\begin{lemma} \label{lem:log-bounds}
For all integers $n \ge 1$ and all $0 \le \alpha \le 1$,
\begin{equation} \label{eq:log_bounds_formula}
    \frac{2\alpha}{2n+\alpha}  \le \log\left(\frac{n+\alpha}{n}\right) \le \frac{\alpha}{n} -\frac{\alpha^2}{2n^2} + \frac{\alpha^3}{3n^3}.
\end{equation}
\end{lemma}

\begin{proof}
\textbf{[Upper bound].} Considering $x = \frac{\alpha}{n}$ and truncating the standard Maclaurin series expansion of $\log(1+x)$ after the first three terms, we immediately obtain:
\begin{equation*}
    \log\left(\frac{n+\alpha}{n}\right) = \log\left(1 + \frac{\alpha}{n}\right) \le \frac{\alpha}{n} -\frac{\alpha^2}{2n^2} + \frac{\alpha^3}{3n^3}.
\end{equation*}
\textbf{[Lower bound].} Considering $x=\frac{\alpha}{2n+\alpha}$ in the classical identity~\cite{Robbins1955ARO,Arslanagic18} for $\log(\frac{1+x}{1-x}) = 2\sum_{k=0}^{\infty} \frac{x^{2k+1}}{2k+1}$, we find:
\begin{equation*}
   \log\left(\frac{n+\alpha}{n}\right) = \frac{2\alpha}{2n+\alpha} + \frac{2}{3}\left(\frac{2\alpha}{2n+\alpha}\right)^3 + \frac{2}{5}\left(\frac{2\alpha}{2n+\alpha}\right)^5 +\ldots  
\end{equation*}
Since every term in the right-hand side is non-negative for $\alpha \in [0,1]$, discarding the higher-order terms preserves the inequality, yielding $\frac{2\alpha}{2n+\alpha} \le \log\left(\frac{n+\alpha}{n}\right)$, which completes the proof.
\end{proof}

\begin{lemma} \label{lem:exp-bound}
For all integers $n \ge 1$ and all $0 \le \alpha \le 1$,
\begin{equation} \label{eq:exp_bounds_formula}
    e^{- \alpha - \frac{\alpha(1-\alpha)}{4n} - \frac{1}{16n}} \;\le\; \left(\frac{n}{n+\alpha}\right)^{n+\frac{1}{2}} \;\le\; e^{-\alpha - \frac{\alpha(1-\alpha)}{4n}}.
\end{equation}
\end{lemma}

\begin{proof}
We first establish bounds for the expression $\left(n+\frac{1}{2}\right)\log\left(\frac{n}{n+\alpha}\right)$. By applying Lemma \ref{lem:log-bounds} and multiplying each term by $-(n+1/2)$, we obtain:
\begin{equation*}
    -\alpha - \frac{\alpha\left(1-\alpha\right)}{2n} - \frac{\alpha^2(4\alpha-3)}{12n^2} - \frac{\alpha^3}{6n^3} \le \left(n+\frac{1}{2}\right)\log\left(\frac{n}{n+\alpha}\right) \le -\frac{(2n+1)\alpha}{2n+\alpha}.
\end{equation*}
In particular, note that $1 \le \left(n+\frac{1}{2}\right)\log\left(\frac{n+1}{n}\right)$ and $\left(n+\frac{1}{2}\right)\log\left(1+\frac{1}{2n}\right) \le \frac{1}{2}+\frac{1}{8n}$. 

Let us define the auxiliary function $g_{n,f}(\alpha)$ over $[0,1]$ as:
\begin{equation*}
    g_{n,f}(\alpha) = \left(n +\frac{1}{2}\right)\log\left(\frac{n+\alpha}{n}\right)-\alpha - \frac{\alpha(1-\alpha)}{4n} - f(n).
\end{equation*}
Differentiating with respect to $\alpha$ yields:
\begin{equation*}
    g'_{n,f}(\alpha)=\frac{n+\frac{1}{2}}{n+\alpha} - 1 - \frac{1}{4n} + \frac{\alpha}{2n}=\frac{(1/2-\alpha)(n-\alpha)}{2n(n+\alpha)}.
\end{equation*}
Since $n-\alpha\geq0$ for all $n \ge 1$, the derivative shows that $g_{n,f}(\alpha)$ is strictly increasing on $[0,1/2]$ and strictly decreasing on $[1/2,1]$. Therefore, it reaches its absolute maximum at $\alpha=1/2$ and its minimum at one of the boundaries $\alpha=0$ or $\alpha=1$. Evaluating $g_{n,f}(\alpha)$ at these critical points gives:
\begin{align*}            
    g_{n,f}(1/2) &= \left(n+\frac{1}{2}\right)\log\left(\frac{2n+1}{2n}\right) - \frac{1}{2}-\frac{1}{16n} - f(n) \le \frac{1}{16n} - f(n), \\
    g_{n,f}(0) &= -f(n), \\
    g_{n,f}(1) &= \left(n +\frac{1}{2}\right)\log\left(1+\frac{1}{n}\right) - 1-f(n) \ge 1 - 1 - f(n) = -f(n).
\end{align*}

\textbf{[Upper bound].} Setting $f(n)=0$, we observe that $g_{n,0}(0) = 0$ and $g_{n,0}(1) \ge 0$. Since the absolute minimum occurs at the boundaries, it follows that $g_{n,0}(\alpha) \ge 0$ for all $0 \le \alpha \le 1$. Hence, $\alpha + \frac{\alpha(1-\alpha)}{4n} \le \left(n+\frac{1}{2}\right)\log\left(\frac{n+\alpha}{n}\right)$, which directly implies:
\begin{equation*}
    \left(n+\frac{1}{2}\right)\log\left(\frac{n}{n+\alpha}\right) \;\le\; -\alpha - \frac{\alpha(1-\alpha)}{4n}.
\end{equation*}

\textbf{[Lower bound].} Setting $f(n) = \frac{1}{16n}$, the absolute maximum satisfies $g_{n,f}(1/2) \le \frac{1}{16n} - \frac{1}{16n} = 0$. Since $\alpha=1/2$ is the maximizer, it follows that $g_{n,f}(\alpha) \le 0$ across the entire interval, meaning:
\begin{equation*}
    \left(n+\frac{1}{2}\right)\log\left(\frac{n+\alpha}{n}\right) \le \alpha + \frac{\alpha(1-\alpha)}{4n} + \frac{1}{16n}.
\end{equation*}
Multiplying by $-1$ isolates the target lower bound:
\begin{equation*}
    - \alpha - \frac{\alpha(1-\alpha)}{4n} - \frac{1}{16n} \;\le\; \left(n+\frac{1}{2}\right)\log\left(\frac{n}{n+\alpha}\right).
\end{equation*}
Taking the exponential function of both uniform bounds yields \eqref{eq:exp_bounds_formula}, completing the proof.
\end{proof}

By synthesizing the uniform control obtained over the localized exponential core in Lemma \ref{lem:exp-bound} with the standard explicit error estimates for Stirling's approximation, we complete the proof for the target factorial ratio.

\begin{proposition}[Stirling Refinement for the Localized Factorial Ratio] \label{prop:n_(n+a)}
For all integers $n\ge 1$ and all real numbers $0\le \alpha\le 1$, 
\begin{equation} \label{eq:stirling_refinement_formula}
    \sqrt{2\pi}\frac{e^{\frac{1}{12n+1} - \frac{\alpha(1-\alpha)}{4n} - \frac{1}{16n}}}{\sqrt{n+\alpha}} < \frac{n!e^{n+\alpha}}{(n+\alpha)^{n+1}} < \sqrt{2\pi}\frac{e^{\frac{1}{12n} - \frac{\alpha(1-\alpha)}{4n}}}{\sqrt{n+\alpha}}.
\end{equation}
\end{proposition}

\begin{proof}
We start from the classical explicit non-asymptotic bounds for Stirling's formula~\cite{Robbins1955ARO} for $n!$:
\begin{equation*}
    \sqrt{2\pi n}\left(\frac{n}{e}\right)^n e^{\frac{1}{12n+1}} < n! < \sqrt{2\pi n}\left(\frac{n}{e}\right)^n e^{\frac{1}{12n}}.
\end{equation*}
Dividing each term by $(n+\alpha)^{n+1}$ yields:
\begin{equation*}
    \sqrt{2\pi n}\left(\frac{n}{e}\right)^n \frac{e^{\frac{1}{12n+1}}}{(n+\alpha)^{n+1}} < \frac{n!}{(n+\alpha)^{n+1}} < \sqrt{2\pi n}\left(\frac{n}{e}\right)^n \frac{e^{\frac{1}{12n}}}{(n+\alpha)^{n+1}}.
\end{equation*}
By grouping the exponential bases and separating the root terms, we can collect the shifting factors as follows:
\begin{equation*}
    \frac{\sqrt{2\pi}}{e^{n}}\frac{e^{\frac{1}{12n+1}}}{\sqrt{n+\alpha}} \left(\frac{n}{n+\alpha}\right)^{n+\frac{1}{2}} < \frac{n!}{(n+\alpha)^{n+1}} <  \frac{\sqrt{2\pi}}{e^{n}}\frac{e^{\frac{1}{12n}}}{\sqrt{n+\alpha}} \left(\frac{n}{n+\alpha}\right)^{n+\frac{1}{2}}.
\end{equation*}
Finally, substituting the sharp lower and upper estimates established in Lemma \ref{lem:exp-bound} into the shifting core leads to:
\begin{equation*}
    \frac{\sqrt{2\pi}}{e^{n}}\frac{e^{\frac{1}{12n+1}}}{\sqrt{n+\alpha}} e^{- \alpha - \frac{\alpha(1-\alpha)}{4n} - \frac{1}{16n}}  < \frac{n!}{(n+\alpha)^{n+1}} < \frac{\sqrt{2\pi}}{e^{n}}\frac{e^{\frac{1}{12n}}}{\sqrt{n+\alpha}} e^{-\alpha - \frac{\alpha(1-\alpha)}{4n}}.
\end{equation*}
Multiplying the entire chain of inequalities by $e^{n+\alpha}$ cancels out the legacy $e^{-n}$ and $e^{-\alpha}$ factors, directly yielding the claimed bounds in \eqref{eq:stirling_refinement_formula}. This holds for all $n \ge 1$ and $0 \le \alpha \le 1$.
\end{proof}

\section{Proofs of the Main Results} \label{sec:proofs}

In this section, we assemble the analytical machinery developed in the preceding sections to complete the proofs of our main theorems. We begin by establishing the sharp explicit bounds for the discrete error sequence $\varepsilon_n$ stated in Theorem \ref{thm:varepsilon_n}. The proof employs a two-pronged approach: a certified interval-arithmetic verification for the initial cases $1 \le n \le 10$, followed by a rigorous inductive argument using proof by contradiction for the infinite tail $n \ge 11$.

\begin{proof}[Proof of Theorem \ref{thm:varepsilon_n}]
The leftmost inequality follows immediately from the definition of $h_n$. To verify the remaining bounds, let us define the auxiliary lower correction sequence $h_n$ for $n \ge 1$ as:
\begin{equation} \label{eq:hn_definition}
    h_n = \frac{\sqrt{2\pi}\left(\sqrt{n+1}-1\right)}{3\sqrt{n}\sqrt{n+1}\left(\sqrt{n}+\sqrt{n+1}\right)}.
\end{equation}

For the initial cases $1 \le n \le 10$, the strict inequalities are rigorously verified via certified numerical computations using the Arb library for interval arithmetic. The values of $\varepsilon_n$ are tightly enclosed in intervals using the convergent series representation of the exponential integral combined with an explicit bound for the truncated tail (see~ \ref{sec:appendix-code}). As showcased in Table \ref{tab:stieltjes-error-bounds}, the resulting verified intervals establish the validity of the bounds for all initial discrete points.

\begin{table}[ht]
\centering
\caption{Numerical comparison of the Stieltjes approximation error with the proposed upper and lower bounds for the initial cases verified via Arb.}
\label{tab:stieltjes-error-bounds}
\begin{tabular}{rccc}
\toprule
$n$ & $\dfrac{\sqrt{2\pi}}{3\sqrt{n}}-h_n$ & $\varepsilon_n$ & $\dfrac{\sqrt{2\pi}}{3\sqrt{n}}$ \\
\midrule
1  & 0.73417442372548448 & 0.82316401210310848 & 0.83554275821033350 \\
2  & 0.51145111449534392 & 0.58755771819609751 & 0.59081795030183868 \\
3  & 0.41777137910516675 & 0.48089027843485560 & 0.48240083637217844 \\
4  & 0.36325423259448364 & 0.41688374523893026 & 0.41777137910516675 \\
5  & 0.32647478488412551 & 0.37307276374895673 & 0.37366608109319525 \\
6  & 0.29946575882540319 & 0.34067926341357864 & 0.34110890264882950 \\
7  & 0.27851195870579887 & 0.31547702692815313 & 0.31580547828364487 \\
8  & 0.26161953032263633 & 0.29514789074208149 & 0.29540897515091934 \\
9  & 0.24761004437548599 & 0.27830053945988535 & 0.27851425273677783 \\
10 & 0.23573613957543276 & 0.26404284438183366 & 0.26422181984040073 \\
\bottomrule
\end{tabular}
\end{table}

For the remaining cases where $n \ge 11$, we proceed analytically and prove each uniform bound independently by induction through contradiction.

\medskip
\noindent\textbf{[Upper bound]} Suppose that there exists an integer $n \ge 10$ which acts as the minimum counterexample index where the upper bound fails at the next step, meaning that:
\begin{equation*}
    \frac{\sqrt{2\pi}}{3\sqrt{n}} \ge \varepsilon_n \qquad \text{and} \qquad \varepsilon_{n+1} > \frac{\sqrt{2\pi}}{3\sqrt{n+1}}.
\end{equation*}
Utilizing the fundamental recurrence identity linking successive error steps, $\varepsilon_n = \varepsilon_{n+1} + \Delta_{n+1}$, we can substitute our assumption to find:
\begin{equation*}
    \frac{\sqrt{2\pi}}{3\sqrt{n}} \ge \varepsilon_{n}  > \frac{\sqrt{2\pi}}{3\sqrt{n+1}} + \Delta_{n+1}.
\end{equation*}
Isolating the endpoint decrement $\Delta_{n+1}$ yields the upper estimate:
\begin{equation} \label{eq:proof_upper_contradiction_1}
    \Delta_{n+1} < \frac{\sqrt{2\pi}}{3\sqrt{n}} - \frac{\sqrt{2\pi}}{3\sqrt{n+1}} < \frac{\sqrt{2\pi}}{3\sqrt{n}\sqrt{n+1}\left(\sqrt{n}+\sqrt{n+1}\right)}.
\end{equation}
Conversely, by Lemma \ref{lem:Delta_n}, the decrement is unconditionally bounded from below by the kernel and residual combination. Reorganizing those terms, we have:
\begin{align*}
     \Delta_{n+1} &\ge \sqrt{2\pi}\left(\frac{\sqrt{n+1}}{3\sqrt{n}(n+1)\left(\sqrt{n+1} + \sqrt{n}\right)} - \frac{1}{n^{5/2}} \left( \frac{7}{1440}+\frac{1}{8n} \right)\right) \\
                  &\ge \frac{\sqrt{2\pi}}{\sqrt{n}}\left(\frac{\sqrt{n+1}}{3(n+1)\left(\sqrt{n+1} + \sqrt{n}\right)} - \frac{1}{n^{2}} \left( \frac{7}{1440}+\frac{1}{8n} \right)\right).
\end{align*}
Since $n \ge 10$, we can explicitly split the numerator of the first term to isolate the baseline factor:
\begin{equation*}
    \Delta_{n+1} \ge \frac{\sqrt{2\pi}}{3\sqrt{n}(n+1)\left(\sqrt{n+1} + \sqrt{n}\right)} + T_n,
\end{equation*}
where the trailing residual remainder $T_n$ is defined as:
\begin{equation*}
    T_n =  \frac{\sqrt{2\pi}}{\sqrt{n}}\left[\frac{2}{3(n+1)\left(\sqrt{n+1} + \sqrt{n}\right)} - \frac{1}{n^{2}} \left( \frac{7}{1440}+\frac{1}{8n} \right)\right].
\end{equation*}
By Lemma \ref{lem:fractions-low}, the expression enclosed within the square brackets is strictly positive for all $n \ge 10$, ensuring that $T_n > 0$. This directly implies:
\begin{equation} \label{eq:proof_upper_contradiction_2}
    \Delta_{n+1} > \frac{\sqrt{2\pi}}{3\sqrt{n}(n+1)\left(\sqrt{n+1} + \sqrt{n}\right)}.
\end{equation}
Equations \eqref{eq:proof_upper_contradiction_1} and \eqref{eq:proof_upper_contradiction_2} stand in direct contradiction on the upper bound capacity of $\Delta_{n+1}$. Consequently, the assumption is false, and $\varepsilon_{n+1} \le \frac{\sqrt{2\pi}}{3\sqrt{n+1}}$ holds for all $n \ge 10$.

\medskip
\noindent\textbf{[Lower bound]} Symmetrically, suppose that there exists a minimum integer $n \ge 10$ such that the proposed lower bound fails at the subsequent step:
\begin{equation*}
    \frac{\sqrt{2\pi}}{3\sqrt{n}} - h_n \le \varepsilon_n \qquad \text{and} \qquad \varepsilon_{n+1} < \frac{\sqrt{2\pi}}{3\sqrt{n+1}} - h_{n+1}.
\end{equation*}
Applying the recurrence relation $\varepsilon_n = \varepsilon_{n+1} + \Delta_{n+1}$ to this joint condition yields:
\begin{equation*}
    \frac{\sqrt{2\pi}}{3\sqrt{n}} - h_n \le \varepsilon_{n}  < \frac{\sqrt{2\pi}}{3\sqrt{n+1}} - h_{n+1} + \Delta_{n+1}.
\end{equation*}
Rearranging the terms to isolate the decrement forces the inequality:
\begin{equation*}
    \Delta_{n+1} > \frac{\sqrt{2\pi}}{3\sqrt{n}} - \frac{\sqrt{2\pi}}{3\sqrt{n+1}} - h_n + h_{n+1} > \frac{\sqrt{2\pi}\sqrt{n+1}}{3\sqrt{n}\sqrt{n+1}\left(\sqrt{n}+\sqrt{n+1}\right)} + h_{n+1},
\end{equation*}
which can be expressed cleanly as:
\begin{equation*}
    \Delta_{n+1} - \frac{\sqrt{2\pi}\sqrt{n+1}}{3\sqrt{n}\sqrt{n+1}\left(\sqrt{n}+\sqrt{n+1}\right)} > h_{n+1}.
\end{equation*}
Invoking the explicit fractional algebraic inequality from Lemma \ref{lem:fraction-high}, we can bound $h_{n+1}$ from below to conclude that:
\begin{equation} \label{eq:proof_lower_contradiction_3}
    \Delta_{n+1} - \frac{\sqrt{2\pi}\sqrt{n+1}}{3\sqrt{n}\sqrt{n+1}\left(\sqrt{n}+\sqrt{n+1}\right)} > \frac{\sqrt{2\pi}}{n^{5/2}}\left(\frac{49}{720}+\frac{97}{1000n}\right).
\end{equation}
However, according to the upper bound established for the decrements in Lemma \ref{lem:Delta_n}, we have:
\begin{equation*}
    \Delta_{n+1} \le \frac{\sqrt{2\pi}\sqrt{n+1}}{3\sqrt{n}(n+1)\left(\sqrt{n+1} + \sqrt{n}\right)} + \frac{\sqrt{2\pi}}{n^{5/2}}\left(\frac{49}{720}+\frac{97}{1000n}\right),
\end{equation*}
which simplifies directly to:
\begin{equation} \label{eq:proof_lower_contradiction_4}
    \Delta_{n+1} - \frac{\sqrt{2\pi}\sqrt{n+1}}{3\sqrt{n}(n+1)\left(\sqrt{n+1} + \sqrt{n}\right)} \le \frac{\sqrt{2\pi}}{n^{5/2}}\left(\frac{49}{720}+\frac{97}{1000n}\right).
\end{equation}
Comparing \eqref{eq:proof_lower_contradiction_3} and \eqref{eq:proof_lower_contradiction_4} reveals a clear structural contradiction. Thus, the lower bound assumption cannot fail, establishing that $\varepsilon_{n+1} \ge \frac{\sqrt{2\pi}}{3\sqrt{n+1}} - h_{n+1}$ for all $n \ge 10$, which completes the proof of the theorem.
\end{proof}

As an immediate consequence of the sharp non-asymptotic bounds proven in Theorem \ref{thm:varepsilon_n}, we formalize the precise rate of convergence of the discrete error.

\begin{corollary} \label{cor:varepsilon_n}
As $n\to\infty$, the error of the Stieltjes approximation satisfies
\[ \varepsilon_n = \frac{\sqrt{2\pi}}{3\sqrt{n}} + O\left(\frac{1}{n}\right). \]
More precisely, for all $n\geq 1$,
\[ 0\leq \frac{\sqrt{2\pi}}{3\sqrt{n}}-\varepsilon_n \leq \frac{\sqrt{2\pi}}{6n}. \]
\end{corollary}
    
\begin{proof}
By Theorem~\ref{thm:varepsilon_n}, we have
\[ 0\leq \frac{\sqrt{2\pi}}{3\sqrt{n}}-\varepsilon_n \leq h_n, \]
where
\[ h_n = \frac{\sqrt{2\pi}(\sqrt{n+1}-1)}{3\sqrt{n}\sqrt{n+1}(\sqrt{n}+\sqrt{n+1})}. \]
Since
\[ \frac{\sqrt{n+1}-1}{\sqrt{n+1}(\sqrt{n}+\sqrt{n+1})} < \frac{1}{\sqrt{n}+\sqrt{n+1}} \leq \frac{1}{2\sqrt{n}}, \]
it follows that
\[ h_n < \frac{\sqrt{2\pi}}{6n}. \]
Consequently,
\[ 0\leq \frac{\sqrt{2\pi}}{3\sqrt{n}}-\varepsilon_n < \frac{\sqrt{2\pi}}{6n}, \]
which yields the asserted asymptotic formula.
\end{proof}

\begin{thm}[Asymptotic Behavior of the Stieltjes--van Boven Approximation Error] \label{thm:sv_main}
For every real number $x \ge e$, the Stieltjes--van Boven approximation error $\varepsilon(x) = \operatorname{li}_{*}(x) - \operatorname{li}(x)$ satisfies the global asymptotic expansion:
\begin{equation} \label{eq:sv_asymptotic_formula}
    \varepsilon(x) = \frac{1}{3} \sqrt{\frac{2\pi}{\log x}} + O\left(\frac{1}{\log x}\right), \qquad x \to \infty.
\end{equation}
\end{thm}

\begin{proof}
We first establish the positivity and strict monotonicity of $\varepsilon(x)$ on its entire domain.

Let $n\geq 1$ and consider $x\in(e^n,e^{n+1})$. By Proposition~\ref{prop:varepsilon(x)-delta(x)}, we have

\[ \varepsilon(x)=\varepsilon_n-\Delta(x).\]

Moreover, Proposition~\ref{prop:Delta(x)order} establishes that $\Delta(x)$ is strictly increasing on $(e^n,e^{n+1})$. Since

\[\Delta(e^{n+1})=\Delta_{n+1},\]

it follows that

\[ \varepsilon(x) > \varepsilon_n-\Delta_{n+1} = \varepsilon_{n+1}. \]

By Theorem~\ref{thm:varepsilon_n}, $\varepsilon_{n+1}>0$, and hence

\[ \varepsilon(x)>0. \]

At the endpoints, positivity follows directly from $\varepsilon(e^n)=\varepsilon_n>0$.

Furthermore, since $\Delta(x)$ is strictly increasing on each interval $(e^n,e^{n+1})$, the function
$\varepsilon(x)$ is strictly decreasing there. At consecutive endpoints, Theorem 1 gives
\[ \varepsilon(e^{n+1}) =\varepsilon_{n+1} < \varepsilon_n  =\varepsilon(e^n). \]

Consequently, combining the strict monotonicity on each interval with the strict decrease of the sequence
$\{\varepsilon_n\}_{n\geq 1}$, we conclude that $\varepsilon(x)$ is strictly decreasing on $[e,\infty)$.

Let $x \ge e$ and set $n=\lfloor\log x\rfloor$. By exploiting the exact decomposition identity from Proposition \ref{prop:varepsilon(x)-delta(x)} alongside the strict monotonicity of the partial error function established in Proposition \ref{prop:Delta(x)order}, the continuous error $\varepsilon(x)$ is tightly squeezed between two successive lattice endpoints:
\begin{equation} \label{eq:sandwich_error}
    \varepsilon_{n+1} = \varepsilon_n - \Delta_{n+1} \le \varepsilon(x) \le \varepsilon_n.
\end{equation}

By direct application of Corollary \ref{cor:varepsilon_n}, the discrete error sequence for every integer $n \ge 1$ evaluates to:
\begin{equation*}
    \varepsilon_n = \frac{\sqrt{2\pi}}{3\sqrt n} + O\left(\frac{1}{n}\right).
\end{equation*}
Since $n=\lfloor\log x\rfloor \sim \log x$ as $x \to \infty$, expanding the reciprocal root yields the smooth baseline transformation:
\begin{equation*}
    \frac{\sqrt{2\pi}}{3\sqrt n} = \frac{1}{3}\sqrt{\frac{2\pi}{\log x}} + O\left(\frac{1}{(\log x)^{3/2}}\right).
\end{equation*}
Symmetrically, evaluating the subsequent step under the same discrete rate gives:
\begin{equation*}
    \varepsilon_{n+1} = \frac{\sqrt{2\pi}}{3\sqrt{n+1}} + O\left(\frac{1}{n+1}\right).
\end{equation*}

Given that $n=\lfloor\log x\rfloor$, we can write the floor index and its shift as localized perturbations of the real domain, namely $n = \log x + O(1)$ and $n+1 = \log x + O(1)$. Applying a standard Taylor expansion to the inverse roots as $x \to \infty$, we obtain the uniform expansions:
\begin{align*}
    \frac{1}{\sqrt n} &= \frac{1}{\sqrt{\log x}} + O\left(\frac{1}{(\log x)^{3/2}}\right), \\
    \frac{1}{\sqrt{n+1}} &= \frac{1}{\sqrt{\log x}}  + O\left(\frac{1}{(\log x)^{3/2}}\right).
\end{align*}
Furthermore, the legacy error terms satisfy the standard big-O asymptotic rates:
\begin{equation*}
    \frac{1}{n}=O\left(\frac{1}{\log x}\right), \qquad \frac{1}{n+1}=O\left(\frac{1}{\log x}\right).
\end{equation*}

Substituting these uniform estimates back into the discrete boundary limits controls both sides of the inequality:
\begin{align*}
    \varepsilon_n &= \frac{\sqrt{2\pi}}{3\sqrt{\log x}} + O\left(\frac{1}{\log x}\right), \\
    \varepsilon_{n+1} &= \frac{\sqrt{2\pi}}{3\sqrt{\log x}} + O\left(\frac{1}{\log x}\right).
\end{align*}
Finally, merging these matching asymptotic behaviors into the squeeze condition \eqref{eq:sandwich_error} forces the continuous error function to converge to the same rate:
\begin{equation*}
    \varepsilon(x) = \frac{1}{3}\sqrt{\frac{2\pi}{\log x}} + O\left(\frac{1}{\log x}\right), \qquad x\to\infty,
\end{equation*}
which completes the formal proof of the theorem.
\end{proof}

\begin{thm}[Asymptotic Behavior of the Absolute Pure Stieltjes Error]
Let $x \ge e$, $n=\lfloor\log x\rfloor$, and $\alpha=\log x-n$. As \(x\to\infty\), the absolute approximation error associated with the pure Stieltjes expansion satisfies
\[|\varepsilon_{S}(x)| = \frac{\sqrt{2\pi}\left|\frac13-\alpha\right|}{\sqrt{\log x}} + O\left(\frac{1}{\log x}\right).\]
\end{thm}

\begin{proof}
By comparing the expansions of the Stieltjes--van Boven and the pure Stieltjes approximations, the relationship between their respective errors can be established via the fractional remainder term:
\[ \varepsilon_{S}(x)=\varepsilon(x)-R(x), \]
where
\[ R(x) = \frac{(\log x-n)x n!}{\log^{n+1}x} = \frac{\alpha n!e^{n+\alpha}}{(n+\alpha)^{n+1}}. \]

From the continuous asymptotic behavior of the Stieltjes--van Boven error established in Theorem \ref{thm:sv_main}, we have
\[ \varepsilon(x) = \frac{\sqrt{2\pi}}{3\sqrt{\log x}} + O\left(\frac{1}{\log x}\right). \]

We now estimate the contribution of the remainder term \(R(x)\). Invoking the sharp uniform factorial ratio estimates from Proposition \ref{prop:n_(n+a)}, it follows that
\[
    \sqrt{2\pi}\, \frac{e^{\frac{1}{12n+1} - \frac{\alpha(1-\alpha)}{4n} - \frac{1}{16n}} }{\sqrt{n+\alpha}}
    \, < \, \frac{n!e^{n+\alpha}}{(n+\alpha)^{n+1}} \, <\, 
    \sqrt{2\pi}\,\frac{e^{\frac{1}{12n}-\frac{\alpha(1-\alpha)}{4n}}}{\sqrt{n+\alpha}}.
\]
Here, the factors \(e^{-\alpha}\) present in the baseline bounds cancel out exactly with the factor \(e^\alpha\) introduced by multiplication by \(e^{n+\alpha}\). Since the fractional part is strictly bounded by \(0\le\alpha<1\), the exponential exponents in these boundaries are \(O(1/n)\) uniformly with respect to \(\alpha\). Consequently, we can write the ratio as a perturbation of the main algebraic term:
\[ \frac{n!e^{n+\alpha}}{(n+\alpha)^{n+1}} = \frac{\sqrt{2\pi}}{\sqrt{n+\alpha}} \left(1+O\left(\frac1n\right)\right), \]
uniformly for \(0\le\alpha<1\). Therefore, the remainder function expands as
\[ R(x) = \frac{\alpha\sqrt{2\pi}}{\sqrt{n+\alpha}} + O\left(\frac{1}{n^{3/2}}\right). \]

Using the identity \(n+\alpha=\log x\) and the basic relation \(n\sim\log x\) as $x \to \infty$, we obtain the targeted smooth approximation for the remainder:
\[ R(x) = \frac{\alpha\sqrt{2\pi}}{\sqrt{\log x}} + O\left(\frac{1}{(\log x)^{3/2}}\right). \]

Substituting the asymptotic expressions for $\varepsilon(x)$ and $R(x)$ back into the structural difference relation yields
\begin{align*}
    \varepsilon_{S}(x) &= \varepsilon(x)-R(x)\\
                     &= \frac{\sqrt{2\pi}}{3\sqrt{\log x}} - \frac{\alpha\sqrt{2\pi}}{\sqrt{\log x}} + O\left(\frac{1}{\log x}\right)\\
                     &= \frac{\sqrt{2\pi}\left(\frac13-\alpha\right)}{\sqrt{\log x}} + O\left(\frac{1}{\log x}\right).
\end{align*}

Finally, utilizing the reverse triangle inequality \(\bigl||u|-|v|\bigr|\leq |u-v|\), we isolate the absolute error oscillation driven by the fractional remainder:
\begin{align*}
    \left||\varepsilon_{S}(x)| - \frac{\sqrt{2\pi}\left|\frac13-\alpha\right|}{\sqrt{\log x}}\right|
      &\leq \left|\varepsilon_{S}(x) - \frac{\sqrt{2\pi}\left(\frac13-\alpha\right)}{\sqrt{\log x}}\right|\\
      &= O\left(\frac{1}{\log x}\right).
\end{align*}
Consequently, we arrive at the explicit asymptotic description
\[
|\varepsilon_{S}(x)| = \frac{\sqrt{2\pi}\left|\frac13-\alpha\right|} {\sqrt{\log x}} +O \left(\frac{1}{\log x}\right),
\]
which completes the proof of the theorem.
\end{proof}

\subsection{Technical Algebraic Support Lemmas} \label{subsec:algebraic_support}

To ensure the structural integrity and unconditional nature of the contradictions established in Section \ref{sec:proofs}, we isolate the underlying algebraic and fractional inequalities. By evaluating these technical relations independently, we avoid breaking the analytical flow of the primary error derivations.

Our first lemma establishes a useful exact rationalization identity for inverse square roots.

\begin{lemma} \label{lem:fractions}
For all positive integers $k \ge 1$, the following exact identity holds:
\begin{equation} \label{eq:fractions_identity}
    \frac{1}{\sqrt{k+1}} = \frac{1}{\sqrt{k}} - \frac{1}{(k+1)\sqrt{k}\left(1 + \sqrt{\frac{k}{k+1}}\right)}.
\end{equation}
\end{lemma}

\begin{proof}
We initiate from the elementary algebraic identity:
\begin{equation*}
    1 = (k+1)\left(1-\frac{k}{k+1}\right).
\end{equation*}
Multiplying both sides by the conjugate root factor yields:
\begin{equation*}
    1+\sqrt{\frac{k}{k+1}} = \frac{1+\sqrt{\frac{k}{k+1}}}{(k+1)\left(1-\frac{k}{k+1}\right)}.
\end{equation*}
Rationalizing the denominator on the right-hand side, we find:
\begin{equation*}
    1 + \sqrt{\frac{k}{k+1}} = \frac{1}{(k+1)\left(1 - \sqrt{\frac{k}{k+1}}\right)}.
\end{equation*}
Reorganizing the terms leads directly to:
\begin{equation*}
    1 - \sqrt{\frac{k}{k+1}} = \frac{1}{(k+1)\left(1 + \sqrt{\frac{k}{k+1}}\right)}.
\end{equation*}
Finally, dividing the entire equation by $\sqrt{k}$ separates the reciprocal root components:
\begin{equation*}
    \frac{1}{\sqrt{k}} - \frac{1}{\sqrt{k+1}} = \frac{1}{(k+1)\sqrt{k}\left(1 + \sqrt{\frac{k}{k+1}}\right)}.
\end{equation*}
Rearranging the terms directly proves \eqref{eq:fractions_identity}.
\end{proof}

Our second lemma justifies the strict positivity of the auxiliary residual remainder sequence $T_n$ introduced during the inductive upper bound proof.

\begin{lemma} \label{lem:fractions-low}

For every integer $n \ge 10$, the following strict fractional inequality holds:

\begin{equation} \label{eq:fractions_low_ineq}
    \frac{2}{3(n+1)(\sqrt{n+1}+\sqrt n)} - \frac{1}{n^2}\left(\frac{7}{1440}+\frac{1}{8n}\right) > 0.
\end{equation}
\end{lemma}

\begin{proof}
For every $n \ge 10$, we exploit the monotonic constraints $\sqrt{n+1} \le 2\sqrt{n}$ and $n+1 \le 2n$. Consequently, the denominator of the first term can be aggressively minorized as:
\begin{equation*}
    3(n+1)(\sqrt{n+1}+\sqrt n) \le 3(2n)(3\sqrt n) = 18n^{3/2}.
\end{equation*}
Taking the reciprocal reverses the inequality, yielding:
\begin{equation} \label{eq:low_inter_1}
    \frac{2}{3(n+1)(\sqrt{n+1}+\sqrt n)} \ge \frac{2}{18n^{3/2}} = \frac{1}{9n^{3/2}}.
\end{equation}
On the other hand, since $n \ge 10$, the trailing inner factor satisfies:
\begin{equation*}
    \frac{7}{1440}+\frac{1}{8n} \le \frac{7}{1440}+\frac{1}{80} = \frac{25}{1440} = \frac{5}{288}.
\end{equation*}
Moreover, evaluating the baseline algebraic rates for $n \ge 10$ reveals that:
\begin{equation*}
    \frac{\sqrt n}{9} \ge \frac{\sqrt{10}}{9} > \frac{5}{288}.
\end{equation*}
Dividing this condition by $n^2$ leads directly to the strict dominant relation:

\begin{equation} \label{eq:low_inter_2}
    \frac{1}{9n^{3/2}} > \frac{5}{288n^2} \ge \frac{1}{n^2}\left(\frac{7}{1440}+\frac{1}{8n}\right).
\end{equation}

Combining the intermediate inequalities \eqref{eq:low_inter_1} and \eqref{eq:low_inter_2} establishes the targeted strict positivity asserted in \eqref{eq:fractions_low_ineq} for all $n \ge 10$.
\end{proof}

Our third lemma secures the strict lower boundary control required to drive the inductive contradiction for the pure Stieltjes error lower bound.

\begin{lemma} \label{lem:fraction-high}
For every integer $n \ge 10$, the shifted correction sequence $h_{n+1}$ satisfies:
\begin{equation} \label{eq:fraction_high_ineq}
    h_{n+1} > \frac{\sqrt{2\pi}}{n^{5/2}} \left(\frac{49}{720}+\frac{97}{1000n}\right),
\end{equation}
where $h_n$ is defined as in \eqref{eq:hn_definition}.
\end{lemma}

\begin{proof}
Let $n \ge 10$. Since $\sqrt{n+2} \ge 2$, it follows that $\sqrt{n+2}-1 \ge \frac{\sqrt{n+2}}{2}$. Additionally, using the basic bounds $\sqrt{n+2} \le 2\sqrt{n+1}$, we find that $\sqrt{n+1}+\sqrt{n+2} \le 3\sqrt{n+1}$. Substituting these minorants into the explicit definition of $h_{n+1}$ yields:
\begin{equation} \label{eq:high_inter_1}
\begin{aligned}
    h_{n+1} &= \frac{\sqrt{2\pi}\left(\sqrt{n+2}-1\right)}{3\sqrt{n+1}\sqrt{n+2} \left(\sqrt{n+1}+\sqrt{n+2}\right)} \\ &\ge \frac{\sqrt{2\pi}\sqrt{n+2}}{6\sqrt{n+1}\sqrt{n+2}\left(\sqrt{n+1}+\sqrt{n+2}\right)} \\ &\ge \frac{\sqrt{2\pi}}{18(n+1)}.
\end{aligned}
\end{equation}
Concurrently, for the right-hand side polynomial terms under $n \ge 10$, we have:
\begin{equation} \label{eq:high_inter_2}
    \frac{49}{720}+\frac{97}{1000n} \le \frac{49}{720}+\frac{97}{10000} < \frac{1}{10}.
\end{equation}
Therefore, to validate \eqref{eq:fraction_high_ineq}, it suffices to demonstrate the strict algebraic inequality:
\begin{equation*}
    \frac{1}{18(n+1)} > \frac{1}{10n^{5/2}} \quad \iff \quad 10n^{5/2} > 18(n+1).
\end{equation*}

This final condition can be rewritten as 

\[10n^{5/2} > 18n\left(1+\frac{1}{n}\right),\]

which simplifies to checking whether (dividing by $10n$ both sides of the inequality)

\[n^{3/2} > \frac{9}{5}\left(1+\frac{1}{n}\right).\]

Since for $n \ge 10$ we have $n^{3/2} \ge 10\sqrt{10} \approx 31.62$ and $\frac{9}{5}\left(1+\frac{1}{10}\right) = \frac{99}{50} = 1.98$, the inequality holds by a wide margin. Combining this with \eqref{eq:high_inter_1} and \eqref{eq:high_inter_2} completes the proof.
\end{proof}

\appendix
\section{ Certified verification of Theorem~\ref{thm:varepsilon_n} for $n = 1,\ldots,10$}
\label{sec:appendix-code}

\begin{lstlisting}[language=Python, caption={Certified verification of Theorem~\ref{thm:varepsilon_n} for $n = 1,\ldots,10$}]

# Certified verification of Theorem 1 for n = 1,...,10
# Run with SageMath.

R = RealBallField(200)
K = 200

gamma = R(euler_gamma)
sqrt_2pi = (2 * R.pi()).sqrt()

def epsilon_ball(n):
    N = R(n)

    # Stieltjes--van Boven approximation at x = exp(n)
    S = R(0)
    for k in range(n):
        S += R(factorial(k)) / (N ** (k + 1))

    approximation = N.exp() * S

    # Finite part of the convergent series for Ei(n)
    partial = R(0)
    for k in range(1, K + 1):
        partial += N**k / (R(k) * R(factorial(k)))

    # Rigorous upper bound for the positive tail
    t = N**(K + 1) / (
        R(K + 1) * R(factorial(K + 1))
    )
    tail_bound = t / (1 - N / R(K + 2))

    # The tail belongs to [0, tail_bound]
    tail = R(tail_bound / 2, tail_bound / 2)

    ei_n = gamma + N.log() + partial + tail

    return approximation - ei_n


def lower_bound(n):
    N = R(n)
    return (
        sqrt_2pi / (3 * N.sqrt())
        - sqrt_2pi * ((N + 1).sqrt() - 1)
        / (
            3 * N.sqrt() * (N + 1).sqrt()
            * (N.sqrt() + (N + 1).sqrt())
        )
    )


def upper_bound(n):
    return sqrt_2pi / (3 * R(n).sqrt())


print(" n       epsilon_n             lower bound"
      "             upper bound")
print("-" * 85)

for n in range(1, 11):
    E = epsilon_ball(n)
    L = lower_bound(n)
    U = upper_bound(n)

    # Certified strict inequalities
    assert L.upper() < E.lower()
    assert E.upper() < U.lower()

    print(n, E, L, U)

print("\nAll ten cases have been certified.")

\end{lstlisting}

\section{Use of AI tools.}

    The author used ChatGPT (OpenAI) and Gemini (Google) as auxiliary tools to assist with language editing, manuscript organization, and the discussion and verification of mathematical arguments. All mathematical results, proofs, and conclusions were critically reviewed and validated by the author, who takes full responsibility for the contents of this paper.

\bibliographystyle{elsarticle-num} 
\bibliography{references}

\end{document}